\documentclass{amsart}

\usepackage{amsmath}
\usepackage{amsfonts}
\usepackage{amssymb}
\usepackage{graphicx}

\newtheorem{theorem}{Theorem}
\numberwithin{theorem}{section} \theoremstyle{plain}

\newtheorem{corollary}[theorem]{Corollary}

\newtheorem{question}{Question}
\newtheorem{proposition}[theorem]{Proposition}

\theoremstyle{definition}
\newtheorem{definition}[theorem]{Definition}
\newtheorem{notation}[theorem]{Notation}

\newtheorem{remark}[theorem]{Remark}

\numberwithin{equation}{section}

\begin{document}

\title[Distribution of the traces of Frobenius]
{Distribution of the traces of Frobenius on elliptic curves over function fields}
\author{Am\'{\i}lcar Pacheco}
\address{Universidade Federal do Rio de Janeiro, Departamen\-to de Matem\'a\-ti\-ca Pura, Rua
Guaiaquil 83, Cachambi, 20785-050 Rio de Janeiro, RJ, Brasil}
\email{amilcar@impa.br}\thanks{This work was partially supported by CNPq research grant
300896/91-3 and Pronex \#41.96.0830.00}

\date{\today}
 \maketitle

\section*{Introduction}

Let $C$ be a smooth irreducible projective curve defined over a finite field $\mathbb{F}_{q}$ of
$q$ elements of characteristic $p>3$ and $K=\mathbb{F}_{q}(C)$ its function field and
$\varphi_{\mathcal{E}}:\mathcal{E}\to C$ the minimal regular model of $\mathbf{E}/K$. For each
$P\in C$ denote $\mathcal{E}_P=\varphi^{-1}_{\mathcal{E}}(P)$. The elliptic curve $E/K$ has good
reduction at $P\in C$ if and only if $\mathcal{E}_P$ is an elliptic curve defined over the residue
field $\kappa_P$ of $P$. This field is a finite extension of $\mathbb{F}_q$ of degree $\deg(P)$.
Let $t(\mathcal{E}_P)=q^{\deg(P)}+1-\#\mathcal{E}_P(\kappa_P)$ be the trace of Frobenius at $P$. By
Hasse-Weil's theorem (cf. \cite[Chapter V, Theorem 2.4]{Silverman}), $t(\mathcal{E}_P)$ is the sum
of the inverses of the zeros of the zeta function of $\mathcal{E}_P$. In particular,
$|t(\mathcal{E}_P)|\le 2q^{\deg(P)}$. Let $C_0\subset C$ be the set of points of $C$ at which $E/K$
has good reduction and $C_0(\mathbb{F}_{q^k})$ the subset of $\mathbb{F}_{q^k}$-rational points of
$C_0$.

\begin{question}Let $B\ge 1$ and $t$ be integers and suppose $|t|\le 2q^{B/2}$. How big is
$\pi(B,t)=\#\{P\in C_0\,|\,\deg(P)\le B\text{ and }t(\mathcal{E}_P)=t\}$?
\end{question}

A similar question was originally posed by Lang-Trotter \cite{LangTrotter} for elliptic curves over
$\mathbb{Q}$ and later extended to elliptic curves over number fields \cite{Murty}.

For each $k\le B$ such that $|t|\le 2q^{k/2}$ we start by estimating $\pi(k,t)'=\#\{P\in
C_0\,|\,\deg(P)=k\text{ and }t(\mathcal{E}_P)=t\}$. Let
$\mathcal{E}_P'=\mathcal{E}_P\times_{\kappa_P}\mathbb{F}_{q^k}$ and $\pi(k,t)''=\#\{P\in
C_0(\mathbb{F}_{q^k})\,|\,t(\mathcal{E}_P')=t\}$. The former set is contained in the latter one so
$\pi(k,t)'\le\pi(k,t)''$ and throughout all this paper we actually estimate $\pi(k,t)''$.

\section{Preliminaries}

Observe first that $E/K$ has to be an ordinary elliptic curve, otherwise $j(E)\in\mathbb{F}_{p^2}$
(cf. \cite[Chapter V, Theorem 3.1]{Silverman}) , but this contradicts the fact that $E/K$ is
non-constant.

Let $j_{\mathcal{E}}:C\to\mathbb{P}^1$ be the $j$-map induced from $\varphi_{\mathcal{E}}$. We say
that $P\in C_0$ is good ordinary, respectively good supersingular, if $\mathcal{E}_P$ is an
ordinary, respectively supersingular, elliptic curve. Since the number of supersingular
$j$-invariants in $\overline{\mathbb{F}}_q$ is finite (cf. \cite[Chapter V, Theorem
4.1]{Silverman}), then the number good of supersingular points $P\in C_0$ is also finite and
bounded by an absolute constant. This does not hold for elliptic curves over $\mathbb{Q}$ (cf.
\cite{Elkies}).

So, we will only concentrate on good ordinary $P\in C_0$. Denote $C_0'=\{P\in C_0\,|\,P\text{ is
ordinary}\}$. Let ${E}/\mathbb{F}_q$ be an elliptic curve and $t({E})=q+1-\#{E}(\mathbb{F}_q)$.
Then ${E}$ is supersingular if and only if $p\nmid t({E})$ (cf. \cite[Ex. 5.10]{Silverman}). Till
the end of this note we assume $p\nmid t$.

\section{Estimate of $\pi(k,t)''$}

\begin{notation}Let $I(t)$ be set of the isogeny classes of elliptic curves
${E}/\mathbb{F}_{q^k}$ defined over such that $\#{E}(\mathbb{F}_{q^k})=q^k+1-t$. Let
$\mathfrak{A}_{k,t}$ be the set of $\mathbb{F}_{q^k}$-isomorphism classes $[{E}]$ of ${E}\in I(t)$
and $N(t)=\#\mathfrak{A}_{k,t}$.
\end{notation}

\begin{definition}Let $\Delta<0$ be an integer such that $\Delta\equiv 0\text{ or }1\pmod 4$,
$B(\Delta)=\{\alpha x^2+\beta xy+\gamma
y^2\,|\,\alpha,\beta,\gamma\in\mathbb{Z},\alpha>0\text{ and
}\beta^2-4\alpha\gamma=\Delta\}$ and $b(\Delta)=\{\alpha x^2+\beta xy+\gamma y^2\in
B(\Delta)\,|\,\gcd(\alpha,\beta,\gamma)=1\}$. The group $\text{SL}_2(\mathbb{Z})$ acts on
$B(\Delta)$ via $(\begin{smallmatrix}\alpha &\beta\\ \gamma
&\delta\end{smallmatrix})f(x,y)=f(\alpha x+\beta y,\gamma x+\delta y)$ preserving
$b(\Delta)$. The sets $b(\Delta)/\text{SL}_2(\mathbb{Z})$ and
$B(\Delta)/\text{SL}_2(\mathbb{Z})$ are finite with cardinality $h(\Delta)$ and
$H(\Delta)$, respectively. The numbers $h(\Delta)$ and $H(\Delta)$ are called the class
number  and the Kronecker class number of $\Delta$,  respectively.
\end{definition}

\begin{proposition}\cite[Proposition 2.2]{Schoof}\label{prop1}
Let $\Delta<0$ be an integer such that $\Delta\equiv 0\text{ or }1\pmod 4$ then
\begin{equation}\label{eqn1}
H(\Delta)=\sum_fh(\Delta/f^2),
\end{equation}
where $f$ runs through all positive divisors
of $\Delta$ such that $\Delta/f^2\in\mathbb{Z}$ and $\Delta/f^2\equiv 0\text{ or }1\pmod
4$.
\end{proposition}

\begin{remark}\label{rem1}
Let $\mathcal{O}$ be an imaginary quadratic order with discriminant $\Delta(\mathcal{O})$
and $h_{\mathcal{O}}$ its class number. It follows from the correspondence between binary
quadratic forms and complex quadratic orders that
$h_{\mathcal{O}}=h(\Delta(\mathcal{O}))$, where $\Delta(\mathcal{O})$ denotes the
discriminant of $\mathcal{O}$ \cite[Chap. 2, Section 7, Theorem 4] {BorSha}.
\end{remark}

\begin{proposition}\cite[Theorem 4.5]{Schoof}
\label{prop2}Let ${E}\in I(t)$ and $\mathcal{O}=\text{End}_{\mathbb{F}_{q^k}}({E})$.
$\#\{[{E}']\in\mathfrak{A}_{k,t}\,|
\,\mathcal{O}=\text{End}_{\mathbb{F}_q}({E}')\}=h_{\mathcal{O}}$.
\end{proposition}

\begin{notation}Denote $\mathcal{O}(t^2-4q^k)$ the imaginary quadratic order
 of discriminant $t^2-4q^k$.
 \end{notation}

\begin{corollary}\label{cor1}\cite[Theorem 4.6]{Schoof}
$N(t)=H(t^2-4q^k)$.
\end{corollary}

\begin{proof}By \cite[Theorem 4.3]{Schoof}, since $p\nmid t$, all
imaginary quadratic orders $\mathcal{O}\supset\mathcal{O}(t^2-4q^k)$ occur as
$\mathbb{F}_{q^k}$-endomorphism ring of elliptic curves in $I(t)$. Hence, the result
follows from Propositions \ref{prop1} and \ref{prop2} and Remark \ref{rem1}.
\end{proof}

\begin{theorem}\label{thm1}
$$\pi(k,t)''\le\deg_s(j_{\mathcal{E}})H(t^2-4q^k),$$
where $\deg_s(j_{\mathcal{E}})$ denotes the separable degree of $j_{\mathcal{E}}$.
\end{theorem}

\begin{proof}Let ${C}_0'(\mathbb{F}_{q^k})$ be the set of
$\mathbb{F}_{q^k}$-rational points of $C_0'$ and $\mathcal{C}_{k,t}=\{P\in
C_0'(\mathbb{F}_{q^k})\,|\,t(\mathcal{E}_P')=t\}$, where
$\mathcal{E}_{P}'=\mathcal{E}_P\times_{\kappa_P}\mathbb{F}_{q^k}$. Define
$\psi:\mathcal{C}_{k,t}\to\mathfrak{A}_{k,t}$ by $\psi(P)=[\mathcal{E}_{P}']$ and let
$j(\mathcal{E}_P)$ be the $j$-invariant of $\mathcal{E}_P$.

We claim that $\psi^{-1}([\mathcal{E}_{P}'])\subset j_{\mathcal{E}}^{-1}(j(\mathcal{E}_{P}))$. In fact, if
$Q\in\psi^{-1}([\mathcal{E}_{P}'])$, then there exists an $\mathbb{F}_{q^k}$-isomorphism between $\mathcal{E}_{Q}'$ and
$\mathcal{E}_{P}'$, in particular $j(\mathcal{E}_{Q})=j(\mathcal{E}_{P})$. Hence,
$\#\psi^{-1}([\mathcal{E}_{P}'])\le\#j_{\mathcal{E}}^{-1}(j(\mathcal{E}_{P})) \le\deg_s(j_{\mathcal{E}})$ and
\begin{equation}\begin{split}\label{eqn0}\pi(k,t)''&=\sum_{[E]\in\psi(\mathcal{C}_{k,t})} \#\psi^{-1}([E])
\\ &\le\deg_s(j_{\mathcal{E}})\#\psi(\mathcal{C}_{k,t})\le\deg_s(j_{\mathcal{E}})
H(t^2-4q^k).\end{split}\end{equation}
\end{proof}

\begin{corollary}$$\pi(B,t)\le\left(\sum_{\substack{k\le B\\ |t|\le
2q^k}}H(t^2-4q^k)\right)\deg_s(j_{\mathcal{E}}).$$
\end{corollary}

\begin{remark}\label{remkey}We would like to compute examples in which we could test whether the
bound of Theorem \ref{thm1} is achieved. One good sort of example comes from modular curves.
However, in this case there is almost no control on $j_{\mathcal{E}}$ in contrast with the $j$-map
$J$ naturally associated to the modular problems. Moreover, if we observe the proof of Theorem
\ref{thm1} closely, we notice that we can replace the regular minimal model by any elliptic surface
$\varphi_{\widetilde{\mathcal{E}}}:\widetilde{\mathcal{E}}\to C$ having $\mathbf{E}/K$ as the
generic fiber, defining the notions of good ordinary (good supersigular) points in terms of the
fibers of $\widetilde{\mathcal{E}}\to C$ being smooth ordinary (supersingular) elliptic curves. In
this set-up it makes sense to consider the trace of Frobenius of the fibers of good ordinary
points. We can also consider elliptic curves $\mathbb{E}\to C_1$ in the sense of \cite[Chapter
2]{KatzMazur} defined over an affine subcurve $C_1\subset C$ with generic fiber $\mathbf{E}/K$ and
compute the number (still denoted $\pi(k,t)''$) of $\mathbb{F}_{q^k}$-rational points $P\in C_1$
corresponding to good ordinary fibers $\mathbb{E}_P$ such that $t(\mathbb{E}_P')=t$. The elliptic
curve comes equipped with a $j$-map $J:C\to\mathbb{P}^1$ and we look for conditions for
$\pi(k,t)''$ to be equal to $\deg(J)H(t^2-4q^k)$.
\end{remark}

\section{Affine models}\label{affine}

Let $X$ be a smooth irreducible projective curve over $\mathbb{F}_q$ and $Y\subset X$ an affine
subcurve. Suppose there exists an elliptic curve $\mathbb{E}\to Y$ with generic fiber
$\mathbf{E}/K$ and a map $J:X\to \mathbb{P}^1$ whose restriction to $Y$ is given by $y\mapsto
j(\mathbb{E}_y)$, where $\mathbb{E}_y$ denotes the fiber of $\mathbb{E}\to Y$ at $y$. Let
$Y'=\{y\in Y\,|\,\mathbb{E}_y\text{ is ordinary}\}$, denote $Y'(\mathbb{F}_{q^k})$ the subset of
$\mathbb{F}_{q^k}$-rational points. Given $y\in Y'(\mathbb{F}_{q^k})$, let $\kappa_y$ be its
residue field and $\mathbb{E}'_y=\mathbb{E}_y\times_{\kappa_y}\mathbb{F}_{q^k}$. Let
$\mathcal{Y}_{k,t}=\{y\in Y'(\mathbb{F}_{q^k})\,|\,t(\mathbb{E}_y')=t\}$ and
$\pi(k,t)''=\#\mathcal{Y}_{k,t}$. Let $\vartheta:\mathcal{Y}_{k,t}\to\mathfrak{A}_{k,t}$ be the map
defined by $y\mapsto [\mathbb{E}_y']$.


\begin{proposition}\label{propkey}
Suppose the following three conditions are satisfied:
\begin{enumerate}
\item $\vartheta^{-1}([\mathbb{E}_{y}'])=J^{-1}(j(\mathbb{E}_{y}))$.
\item $\vartheta$ is surjective.
\item For every $y\in Y$, the inertia degree $f(y\,|\,j(\mathbb{E}_y))$ equals $1$.
The set $\mathcal{R}\subset Y$ of possible ramification points of $J$ is contained in
$J^{-1}(\{0,1728\})$. For each $y\in\mathcal{R}$ the ramification index $e(y\,|\,0)$ (respectively
$e(y\,|\,1728)$) of $P$ over $0$, respectively $1728$, equals $3$, respectively $2$.
\end{enumerate}
Then $\pi(k,t)''=\deg_s(J)H(t^2-4q^k)$, where $\deg_s(J)$ denotes the separable degree of $J$.
\end{proposition}

\begin{proof} In the definition of $H(\Delta)$, we count the
forms $\alpha x^2+\alpha y^2$, respectively $\alpha x^2+\alpha xy+\alpha y^2$, in $B(\Delta)$, if
they occur, with multiplicity $1/2$, respectively $1/3$. Then we need to replace $h(\Delta)$ in
Proposition \ref{prop1} by $h_w(\Delta)$, where $h_w(-3)=1/3$, $h_w(-4)=1/2$ and
$h_w(\Delta)=h(\Delta)$, for $\Delta<-4$. The equality (\ref{eqn1}) does not change when
reinterpreted with these multiplicities \cite[Proposition 2.1]{SchoofVlugt}. Hence, by Propositions
\ref{prop1} and \ref{prop2} and Remark \ref{rem1}, (cf. (\ref{eqn0}))
\begin{equation*}\begin{split}\pi(k,t)''&=\sum_{[{E}]\in\mathfrak{A}_{k,t}}
\#J^{-1}(j({E}))\\
&=\deg_s(J)\sum_{\mathcal{O}(t^2-4q^k)\subset\mathcal{O}} h_w(\Delta(\mathcal{O}))=\deg_s(J)H(t^2-4q^k).
\end{split}\end{equation*}
\end{proof}

\section{Universal elliptic curves}

\subsection{Igusa curves}

Let ${E}$ be an elliptic curve defined over a field $L$ of characteristic $p$. The absolute
Frobenius $F_{\text{abs}}$ induces an isogeny $F_{\text{abs}}:{E}\to{E}^{(p)}$, where ${E}^{(p)}$
denotes the elliptic curve obtained by raising the coefficients of a Weierstrass equation of ${E}$
to the $p$-th power. For each $n\ge 1$, let $F^n_{\text{abs}}:{E}\to{E}^{(p^n)}$ be the $n$-th
iterate of $F^n_{\text{abs}}$. Let $V^n$ be the dual isogeny of the $n$-th iterate
$F^n_{\text{abs}}$ of $F_{\text{abs}}$. An Igusa structure of level $p^n$ in $\mathbf{E}$ is a
generator of $\ker(V^n)$.

There exists a smooth affine curve $Y_n$ over $\mathbb{F}_p$ parametrizing isomorphism classes of
pairs $({E},P)$, where ${E}$ is an elliptic curve defined over an $\mathbb{F}_p$-scheme $S$ and
$P\in {E}^{(p^n)}(S)$ is an Igusa structure  of level $p^n$. In fact, $Y_n$ is a coarse moduli
scheme for the moduli problem $[\text{Ig}(p^n)]:{E}/S/\mathbb{F}_p\mapsto P$. The compactification
$X_n$ of $Y_n$ obtained by adding $\phi(p^n)/2$ points at infinity (called the cusps) is a smooth
projective irreducible curve over $\mathbb{F}_p$ called the Igusa curve of level $p^n$
\cite[Chapter 12]{KatzMazur}, where $\phi$ denotes the Euler function.

An elliptic curve is $E/S/\mathbb{F}_p$ is ordinary, if each of its geometric fibers is ordinary.
An Igusa ordinary (respectively Igusa supersingular) point $y\in Y_n$ is a point representing the
isomorphism of a pair $({E},P)$, where ${E}/S/\mathbb{F}_p$ is an ordinary elliptic curve, $S$ an
$\mathbb{F}_p$-scheme (respectively ${E}/L$ is a supersingular elliptic curve, $L$ a field of
characteristic $p$) and $P\in{E}^{(p^n)}(S)$ (respectively $P\in{E}^{(p^n)}(L)$) is an Igusa
structure of level $p^n$ in ${E}$.

The group $(\mathbb{Z}/p^n\mathbb{Z})^*$ acts on $Y_n$ by $a\mapsto({E},aP)$ and the group $\{\pm
1\}$ acts trivially. These actions are extended to $X_n$ permuting the cusps simply transitively.
Let $y\in Y_n$ represent the isomorphism class of a pair $({E},P)$. If $y$ is Igusa supersingular,
then $y$ is fixed by $(\mathbb{Z}/p^n\mathbb{Z})^*$. If $y$ is Igusa ordinary and $j({E})=1728$,
respectively $j({E})=0$, then $y$ has a stabilizer of order $2$, respectively $3$, in
$(\mathbb{Z}/p^n\mathbb{Z})^*/\{\pm 1\}$. On all other points of $Y_n$,
$(\mathbb{Z}/p^n\mathbb{Z})^*/\{\pm 1\}$ acts freely. We identify the quotient of $X_n$ by
$(\mathbb{Z}/p^n\mathbb{Z})^*/\{\pm 1\}$ to the projective line $\mathbb{P}^1$ and the quotient map
$J_n:X_n\to\mathbb{P}^1$ is Galois of degree $\phi(p^n)/2$. Its restriction to $Y_n$ is given by
$({E},P)\mapsto j({E})$.

The curve $Y_n^{\text{ord}}$ obtained from $Y_n$ by removing the Igusa supersingular points is a
fine moduli space for the restriction of $[\text{Ig}(p^n)]$ to ordinary elliptic curves. This means
that there exists a universal elliptic curve $\mathbb{E}_n\to Y^{\text{ord}}_n$ such that every
ordinary elliptic curve ${E}/S/\mathbb{F}_p$ with an Igusa structure $P\in {E}^{(p^n)}(S)$ of level
$p^n$ is obtained from $\mathbb{E}_n\to Y^{\text{ord}}_n$ by a unique base extension. In
particular, if $K_n$ is the function field of $X_n$ over $\mathbb{F}_p$ and $\mathbf{E}_n/K_n$ is
the generic fiber of $\mathbb{E}_n\to Y^{\text{ord}}_n$, then $\mathbf{E}_n/K_n$ is the unique
elliptic curve defined over $K_n$ with $j$-invariant $j(\mathbf{E}_n)$ and a $K_n$-rational Igusa
structure of level $p^n$.

If ${E}\in I(t)$ and $P\in{E}^{(p^n)}(\mathbb{F}_{p^k})$ is an Igusa structure of level $p^n$.
Since ${E}$ and ${E}^{(p^n)}$ are isogenous, then $t\equiv p^k+1\pmod{p^n}$. So for the rest of
this subsection, we assume $t\equiv p^k+1\pmod{p^n}$.

\begin{proposition}\label{propa}
Conditions (1), (2) and (3) of Proposition \ref{propkey} are satisfied, a fortiori
$\pi(k,t)''=(\phi(p^n)/2)H(t^2-4p^k)$.
\end{proposition}

\begin{proof}In the notation of Section \ref{affine}, $Y=Y_n^{\text{ord}}$.
Let $y\in\mathcal{Y}_{k,t}$, denote $\mathbb{E}_{n,y}/\kappa_y$ the fiber of $\mathbb{E}_n\to
Y_n^{\text{ord}}$ at $y$ and $\mathbb{E}_{n,y}'=\mathbb{E}_{n,y}\times_{\kappa_y}\mathbb{F}_{p^k}$.

(1) Let $x\in\vartheta^{-1}([\mathbb{E}_{n,y}'])$, then $\mathbb{E}_{n,x}'$ is
$\mathbb{F}_{p^k}$-isomorphic to $\mathbb{E}_{n,y}'$, in particular
$j(\mathbb{E}_{n,x})=j(\mathbb{E}_{n,y})$, i.e., $x\in J_n^{-1}(j(\mathbb{E}_{n,y}))$. Let $x\in
J_n^{-1}(j(\mathbb{E}_{n,y}))$, then $x$ represents the isomorphism class of the pair
$(\mathbb{E}_{n,y},P_y)$, where $P_y\in\mathbb{E}_{n,y}^{(p^n)}(\kappa_y)$ is an Igusa structure of
level $p^n$. By the geometric description of $J_n$, there is no inertia at $Y_n^{\text{ord}}$,
hence $\kappa_x=\kappa_y$. Furthermore, $\mathbb{E}_{n,y}$ is an elliptic curve over
$\kappa_y=\kappa_x$ with $j$-invariant equal to $j(\mathbb{E}_{n,y})=j(\mathbb{E}_{n,x})$ having a
$\kappa_y$-rational Igusa structure. It follows from the universal property of $\mathbb{E}_n\to
Y_n^{\text{ord}}$ that $\mathbb{E}_{n,y}=\mathbb{E}_{n,x}$, a fortiori
$[\mathbb{E}_{n,x}']=[\mathbb{E}_{n,y}']$ and $x\in\vartheta^{-1}([\mathbb{E}_{n,y}'])$.

(2) For every $[{E}]\in\mathfrak{A}_{k,t}$,
$\#{E}(\mathbb{F}_{p^k})=\#{E}^{(p^n)}(\mathbb{F}_{p^k})\equiv 0\pmod{p^n}$. Thus, there exists an
Igusa $P\in{E}^{(p^n)}(\mathbb{F}_{p^k})$ structure of level $p^n$. Let $y\in
Y_n^{\text{ord}}(\mathbb{F}_{p^k})$ represent the isomorphism class of the pair $({E},P)$. But
$\mathbb{E}_{n,y}$ is the unique elliptic curve over $\kappa_y$ with $j$-invariant
$j(\mathbb{E}_{n,y})=j(E)$ having a $\kappa_y$-rational Igusa structure of level $p^n$. Thus,
$\mathbb{E}_{n,y}'={E}$. In particular, $[\mathbb{E}_{n,y}']=[{E}]$ and $\vartheta$ is surjective.

Condition (3) follows from the geometric description of $J_n$. Consequently,
$$\pi(k,t)''=\frac{\phi(p^n)}2\sum_{\mathcal{O}(t^2-4p^k)\subset\,\mathcal{O}}
h_w(\Delta(\mathcal{O}))=\frac{\phi(p^n)}2H(t^2-4p^k).$$\end{proof}

\begin{remark}
Proposition \ref{propa} was implicitly used in \cite[Corollary 2.13]{Pacheco} to
obtain an explicit expression for $\#X_n(\mathbb{F}_{p^k})$.
\end{remark}

\subsection{The modular curve $X(N)$}

Let $N>2$ be an integer not divisible by $p$. Let $\zeta\in\overline{\mathbb{F}}_p$ be a primitive
$N$-th root of unity and $\mathbb{F}_q=\mathbb{F}_p(\zeta)$. Let $Y(N)$ be the affine smooth curve
defined over $\mathbb{F}_q$ parametrizing isomorphism classes of triples $({E},P,Q)$, where ${E}$
is an elliptic curve defined over an $\mathbb{F}_q$-scheme $S$ and $P,Q\in{E}[N](S)$ is a Drinfeld
basis for ${E}[N](S)$ and $e_N(P,Q)=\zeta$ \cite[3.1]{KatzMazur}, where $e_N$ denotes the $N$-th
Weil pairing (cf. \cite[2.8]{KatzMazur} and \cite[III, \S8]{Silverman}). In fact, it is a fine
moduli space for the modular problem $[\Gamma(N)]:{E}/S/\mathbb{F}_q\mapsto (P,Q)$ such that
$e_N(P,Q)=\zeta$. The compactification $X(N)$ of $Y(N)$ obtained by adding the cusps is a smooth
projective irreducible curve defined over $\mathbb{F}_q$ \cite[Theorem 3.7.1]{KatzMazur}.

The group $\text{SL}_2(\mathbb{Z}/N\mathbb{Z})$ acts on $Y(N)$ by $(\begin{smallmatrix}a &b\\ c&
d\end{smallmatrix})(E,P,Q)\mapsto ({E},aP+bQ,cP+dQ)$ and the group $\{\pm 1\}$ acts trivially. If
$y\in Y(N)$ represents the isomorphism class of a triple $({E},P,Q)$, then $y$ has a stabilizer of
order $3$, respectively $2$, if $j(E)=0$, respectively $j({E})=1728$. On all other points of
$Y(N)$, $\text{SL}_2(\mathbb{Z}/N\mathbb{Z})/\{\pm 1\}$ acts freely. The stabilizer at every cusp
has order $N$ \cite[Theorem 6]{Igusa}. So we identify the quotient of $X(N)$ by
$\text{SL}_2(\mathbb{Z}/N\mathbb{Z})/\{\pm 1\}$ to the projective line $\mathbb{P}^1$. Let
$\mathbb{E}(N)\to Y(N)$ be the universal elliptic curve of $Y(N)$ and $\mathbf{E}(N)/K(N)$ its
generic fiber. The quotient map $J(N):X(N)\to\mathbb{P}^1$ is Galois of degree
$(1/2)N\phi(N)\psi(N)$ and its restriction to $Y(N)$ is given by $({E},P,Q)\mapsto j({E})$, where
$\psi(N)=N\prod_{\ell\mid N}(1+(1/\ell))$ and $\ell$ runs over the divisors of $N$.

Denote $\mathcal{O}((t^2-4q^k)/N^2)$ the imaginary quadratic order of discriminant
$(t^2-4q^k)/N^2$. Let ${E}\in I(t)$, we have ${E}[N]\subset{E}(\mathbb{F}_{q^k})$ if and only if
$t\equiv q^k+1\pmod{N^2}$, $q^k\equiv 1\pmod N$ and
$\mathcal{O}((t^2-4q^k)/N^2)\subset\text{End}_{\mathbb{F}_{q^k}}({E})$ \cite[Proposition
3.7]{Schoof}.  Assume till the end of this subsection that $t\equiv q^k+1\pmod{N^2}$ and $q^k\equiv
1\pmod N$. Let $\mathfrak{A}_{k,t}'=\{[{E}]\in\mathfrak{A}_{k,t}\,|\,
\mathcal{O}((t^2-4q^k)/N^2)\subset\text{End}_{\mathbb{F}_{p^k}}({E})\}$. By \cite[Theorem
4.9]{Schoof}, $\#\mathfrak{A}_{k,t}'=H((t^2-4q^k)/N^2)$. Note that $H((t^2-4q^k)/N^2)<H(t^2-4q^k)$.

\begin{proposition}Condition (1) and (3) of Proposition \ref{propkey} are satisfied.
However, $\pi(k,t)''=(1/2)N\phi(N)\psi(N)H((t^2-4q^k)/N^2)$.
\end{proposition}

\begin{proof}In the notation of Section \ref{affine}, $Y=Y(N)$.
Let $y\in\mathcal{Y}_{k,t}$, denote $\mathbb{E}(N)_y/\kappa_y$ the fiber of $\mathbb{E}(N)\to Y(N)$
at $y$ and $\mathbb{E}(N)_{y}'=\mathbb{E}(N)_y\times_{\kappa_y}\mathbb{F}_{q^k}$.

(1) Let $x\in\vartheta^{-1}([\mathbb{E}(N)_{y}'])$, then $\mathbb{E}(N)_x'$ is
$\mathbb{F}_{q^k}$-isomorphic to $\mathbb{E}(N)_y'$, in particular
$j(\mathbb{E}(N)_x)=j(\mathbb{E}(N)_y)$, i.e., $x\in J(N)^{-1}(j(\mathbb{E}(N)_y))$. Let $x\in
J(N)^{-1}(j(\mathbb{E}(N)_y))$, then $x$ represents the isomorphism class of the triple
$(\mathbb{E}(N)_y,P_y,Q_y)$, where $P_y,Q_y$ is a basis for $\mathbb{E}(N)_y(\kappa_y)$ with
$e_N(P,Q)=\zeta$. By the geometric description of $J(N)$, there is no inertia, hence
$\kappa_x=\kappa_y$. Furthermore, $\mathbb{E}(N)_y$ is an elliptic curve over $\kappa_y=\kappa_x$
with $j$-invariant $j(\mathbb{E}(N)_y)=j(\mathbb{E}(N)_x)$ and a $\kappa_y$-rational basis
$P_y,Q_y$ of $\mathbb{E}(N)_y[N]$ such that $e_N(P,Q)=\zeta$. By the universal property of
$\mathbb{E}(N)\to Y(N)$, $\mathbb{E}(N)_x=\mathbb{E}(N)_y$, a fortiori
$[\mathbb{E}(N)_x']=[\mathbb{E}(N)_y']$ and $x\in\vartheta^{-1}([\mathbb{E}(N)_y'])$.

 (2) For every ${E}\in\mathfrak{A}_{k,t}'$, by hypothesis
$\mathcal{O}((t^2-4q^k)/N^2)\subset\text{End}_{\mathbb{F}_{q^k}}({E})$, hence ${E}[N]\subset{E} (\mathbb{F}_{q^k})$, in
particular there exists a basis $P,Q$ for ${E}[N](\mathbb{F}_{q^k})$ such that $e_N(P,Q)=\zeta$. Let $y\in
Y(N)(\mathbb{F}_{q^k})$ represent the isomorphism class of the triple $({E},P,Q)$ satisfying $e_N(P,Q)=\zeta$. But,
$\mathbb{E}(N)_y$ is the unique elliptic curve defined over $\kappa_y$ with $\kappa_y$-rational basis $(P_y,Q_y)$ of
$\mathbb{E}(N)_y[N]$ satisfying $e_N(P_y,Q_y)=\zeta$. Thus $\mathbb{E}(N)_{y}'={E}$. In particular,
$[\mathbb{E}(N)_{y}']=[{E}]$ and $\vartheta$ is onto $\mathfrak{A}_{k,t}'$.

Condition (3) follows from the geometric description of $J(N)$. Therefore,


\begin{equation*}\begin{aligned}
\pi(k,t)''&=\frac 12N\phi(N)\psi(N)\sum_{\mathcal{O}((t^2-4q^k)/N^2)\subset\mathcal{O}}
h_w(\Delta(\mathcal{O}))\\&=\frac 12N\phi(N)\psi(N)H(\frac{t^2-4q^k}{N^2}).
\end{aligned}\end{equation*}
\end{proof}

\subsection{The modular curve $X_1(N)$}

Let $N>4$ be an integer not divisible by $2$, $3$ and $p$. Let $Y_1(N)$ be the smooth affine curve
defined over $\mathbb{F}_p$ parametrizing isomorphism classes of pairs $({E},P)$, where ${E}$ is an
elliptic curve $E$ defined over an $\mathbb{F}_p$-scheme $S$ and $P\in E(S)$ is a point  of exact
order $N$ \cite[Chapter 3]{KatzMazur}. In fact, it is a fine moduli space for the moduli problem
$[\Gamma_1(N)]$ defined by $({E}/S/\mathbb{F}_p,P)\mapsto P$. The compactification $X_1(N)$ of
$Y_1(N)$ is a smooth irreducible projective curve defined over $\mathbb{F}_p$ \cite[Theorem
3.7.1]{KatzMazur}.

Let $\mathbb{E}_1(N)\to Y_1(N)$ be the universal elliptic curve of $Y_1(N)$ and
$\mathbf{E}_1(N)/K_1(N)$ its generic fiber. The $j$-map $J(N):X(N)\to\mathbb{P}^1$ factors through
the Galois cover $X(N)\to X_1(N)$ of degree $N$, whose restriction to $Y(N)$ maps to $Y_1(N)$ by
$(E,P,Q)\mapsto (E,P)$. It induces the $j$-map $J_1(N):X_1(N)\to\mathbb{P}^1$ whose restriction to
$Y_1(N)$ is given by $(E,P)\mapsto j(E)$. Since $2,3\nmid N$, given $y\in Y(N)$ and $y_1\in Y_1(N)$
such that $J(N)(y)=J_1(N)(y_1)$ equals $0$, respectively $1728$, then the ramification index
$e(y\,|\,y_1)$ equals $1$. A fortiori, $e(y_1\,|\,0)=3$, respectively $e(y_1\,|\,1728)=2$. Note
also that since there exists no inertia in $Y(N)\to\mathbb{A}^1$, then the same holds for
$Y_1(N)\to\mathbb{A}^1$, thus Condition (3) of Proposition \ref{propkey} is satisfied.

Observe that if ${E}\in I(t)$ has a point $P\in{E}(\mathbb{F}_q)$ of exact order $N$, then
$N\mid\#{E}(\mathbb{F}_{p^k})$. The converse holds if $N$ is a prime number. We assume till the end
of this subsection that $N=\ell$ is a prime number different from $2$, $3$ and $p$, and $t\equiv
p^k+1\pmod{\ell}$.

\begin{proposition}Conditions (1), (2) and (3) of Proposition \ref{propkey} are
satisfied, a fortiori $\pi(k,t)''=(1/2)(\ell^2-1)H(t^2-4p^k)$.
\end{proposition}

\begin{proof}In the notation of Section \ref{affine}, $Y=Y_1(\ell)$.
Let $y\in\mathcal{Y}_{k,t}$, $\mathbb{E}_1(\ell)_y/\kappa_y$ the fiber of
$\mathbb{E}_1(\ell)\to Y_1(\ell)$ at $y$ and
$\mathbb{E}_1(\ell)_y'=\mathbb{E}_1(\ell)_y\times_{\kappa_y}\mathbb{F}_{p^k}$.

(1) Let $x\in\vartheta^{-1}([\mathbb{E}_1(\ell)_y'])$, then $\mathbb{E}_1(\ell)_x'$ is $\mathbb{F}_{p^k}$-isomorphic to
$\mathbb{E}_1(\ell)_y'$, in particular $j(\mathbb{E}_1(\ell)_x)=j(\mathbb{E}_1(\ell)_y)$, i.e., $x\in
J_1(\ell)^{-1}(j(\mathbb{E}_1(\ell)_y))$. Let $x\in J_1(\ell)^{-1}(j(\mathbb{E}_1(\ell)_y))$, then $x$ represent the
isomorphism class of the pair $(\mathbb{E}_1(\ell)_y,P_y)$, where $P_y\in\mathbb{E}_1(\ell)_y($ \linebreak $\kappa_y)$
is a point of exact order $\ell$. By the geometric description of $J_1(\ell)$, there is no inertia, so
$\kappa_x=\kappa_y$. Furthermore, $\mathbb{E}_1(\ell)_y$ is an elliptic curve over $\kappa_y=\kappa_x$ with
$j$-invariant $j(\mathbb{E}_1(\ell)_y)=j(\mathbb{E}_1(\ell)_x)$ and a $\kappa_y$-rational point $P_y$ of exact order
$\ell$. By the universal property of $\mathbb{E}_1(\ell)\to Y_1(\ell)$, $\mathbb{E}_1(\ell)_x=\mathbb{E}_1(\ell)_y$, a
fortiori $[\mathbb{E}_1(\ell)_x']=[\mathbb{E}_1(\ell)_y']$ and $x\in\vartheta^{-1}([\mathbb{E}_1(\ell)_y'])$.

(2) For every ${E}\in\mathfrak{A}_{k,t}$, by hypothesis, $\ell\mid\#{E}(\mathbb{F}_{p^k})$, thus
there exists $P\in{E}(\mathbb{F}_{p^k})$ of exact order $\ell$. Let $y\in
Y_1(\ell)(\mathbb{F}_{p^k})$ represent the isomorphism class of the pair $({E},P)$. But
$\mathbb{E}_1(\ell)_y$ is the unique elliptic curve defined over $\kappa_y$ with a
$\kappa_y$-rational point $P_y$ of exact order $\ell$. Thus $\mathbb{E}_1(\ell)_y'=E$. In
particular, $[\mathbb{E}_1(\ell)_{y}']=[{E}]$ and $\vartheta$ is surjective.

Condition (3) follows from the geometric description of $J_1(\ell)$. Therefore,
$$\pi(k,t)''=\frac
12\phi(\ell)\psi(\ell)\sum_{\mathcal{O}(t^2-4p^k)\subset\,\mathcal{O}}
h_w(\Delta(\mathcal{O}))=\frac 12(\ell^2-1)H(t^2-4p^k).$$
\end{proof}

\end{document}